\documentclass[12pt]{amsart}
\usepackage{latexsym}
\usepackage{amsfonts}
\usepackage{amsthm}
\usepackage{amsmath}
\usepackage{amssymb}

\def\<{{\langle}}
\def\>{{\rangle}}

\def\eps{\epsilon}

\def\note#1{{}}

\def\note#1{}

\def\rend#1#2{{{\rm End}\sb{#1}(#2)}}

\def\beq{\begin{equation}}
\def\eeq{\end{equation}}

\def\id{{I}}

\def\ot{{\otimes}}

\newcounter{zlist}

\def\Label#1{\label{#1}\ifmmode\llap{[#1] }\else
\marginpar{\smash{\hbox{\tiny [#1]}}}\fi}
\def\Label{\label}

\newtheorem{proposition}{Proposition}[section]

\newtheorem{theorem}[proposition]{Theorem}

\theoremstyle{definition}
\newtheorem{definition}[proposition]{Definition}

\theoremstyle{remark}
\newtheorem{remark}[proposition]{Remark}

\newcounter{c}

\newcommand{\etyk}[1]{\vspace{-7.4mm}$$\begin{equation}\Label{#1}
\addtocounter{c}{1}}
\renewcommand{\]}{\ifnum \value{c}=1 $$\else \end{equation}\fi}
\begin{document}

\title{Spectral--parameter dependent Yang--Baxter operators and
Yang--Baxter systems from algebra structures}
\author{Florin F. Nichita}
\address{Institute of Mathematics "Simion Stoilow" of the Romanian Academy, P.O. Box
1-764, RO-014700 Bucharest, Romania}
\email{Florin.Nichita@imar.ro}
\author{Deepak Parashar}
\address{Mathematics Institute, University of Warwick, Coventry CV4 7AL, U.K.}
\address{Department of Mathematics, University of Wales Swansea,
Singleton Park, Swansea SA2 8PP, U.K.}
\email{parashar@maths.warwick.ac.uk}

\begin{abstract}
For any algebra two families of coloured Yang-Baxter operators are constructed,
thus producing solutions to the two-parameter quantum Yang-Baxter equation.
An open problem about a system of functional equations is stated. The matrix
forms of these operators for two and three dimensional
algebras are computed. A FRT bialgebra for one of these families is
presented. Solutions for the one-parameter quantum Yang-Baxter
equation are derived and a Yang-Baxter system constructed.

\noindent \it{Comm. Algebra 34 (2006) 2713--2726.}
\end{abstract}

\maketitle

\section{Introduction}

\smallskip

The quantum Yang-Baxter equation (QYBE) plays a crucial role in analysis
of integrable systems, in quantum and statistical mechanics and also in
the
theory of quantum groups. In the quantum group theory, solutions of
the constant QYBE lead to examples of bialgebras
via the Faddeev--Reshetikhin--Takhtajan (FRT) construction \cite{frt, Kas:Qua}. On the other hand, the
theory of integrable Hamiltonian
systems makes great use of the solutions of the one-parameter form of the
QYBE, since coefficients of the power series expansion of
such a solution give rise to commuting integrals of motion. The
purpose of this paper is two-fold. First
we study non-additive solutions of the two-parameter form of the
QYBE, also known as the {\em coloured} QYBE.
Such a solution is referred to as a {\em coloured
Yang-Baxter operator}. Secondly, we construct the one-parameter Yang-Baxter
operators and
associate to it a Yang-Baxter system. We also compare our results with other well-known solutions.
It is imperative to note
that the Yang-Baxter operators presented here are obtained from
{\em algebra} structures, and are therefore distinct from the
$R$--matrices that arise from quasitriangular Hopf algebras. See
\cite{DasNic:yan,Kas:Qua,LamRad:Int,Nic:sel}
for ordinary Yang-Baxter operators, \cite{LamRad:Int} for
one-parameter form of the QYBE, \cite{HlaWan,LamRad:Int} for the coloured QYBE,
\cite{GonVes:yan, Ves:yan} for
Yang-Baxter maps and \cite{deeps, quesne} for coloured quantum groups.

\bigskip

\section{Construction of coloured Yang--Baxter operators}

\smallskip

Formally, a coloured Yang-Baxter operator is defined as a function $$ R
:X\times X \to \rend k {V\otimes_k V}, $$ where $X$ is a set and $V$ is a
finite dimensional vector space over a field $k$. Thus for any $u,v\in X$,
$R(u,v) : V\otimes_k V\to V\otimes_k V$ is a linear operator. Starting
with this operator one constructs three operators acting on a triple
tensor product $V\otimes_k V\otimes_k V$, $R_{12}(u,v) = R(u,v)\otimes_k
\id$, $R_{23}(v,w)= \id\otimes_k R(v,w)$, and similarly $R_{13}(u,w)$ as
an operator that acts non-trivially on the first and third factor in
$V\otimes_k V\otimes_k V$. $R$ is a coloured Yang-Baxter operator if it
satisfies the two-parameter form of the QYBE,
\begin{equation}\label{yb}
R_{12}(u,v)R_{13}(u,w)R_{23}(v,w) = R_{23}(v,w)
R_{13}(u,w)R_{12}(u,v)
\end{equation}
for all $u,v,w\in X$. This is essentially the spectral parameter dependent
QYBE and non-additivity implies $R(u,v)\neq R(u-v)$. One of the aims of this paper
is to construct two classes of coloured Yang-Baxter operators based on the
use of associative algebra structures built on $V$. We also state an
open problem related to a system of functional equations.
All tensor products appearing in this paper are defined over $k$.
From now on we assume that $X$ is equal to (a subset of) the ground field $k$.
The method of constructing solutions to equation~(\ref{yb})  is based on the
ideas applied in \cite{DasNic:yan,Nic:sel} to the case
of the constant QYBE:
\begin{equation}\label{constyb}
R_{12} R_{13} R_{23} = R_{23} R_{13} R_{12}
\end{equation}
The key point of the construction is to suppose
that $V=A$ is an associative
$k$-algebra, and then to derive a solution to equation~(\ref{yb}) from
the associativity of the product in $A$.
Guided by the observation that the operator
$$
R:A\otimes A\to A\otimes A \quad \text{given by} \quad
a\otimes b\mapsto 1\otimes ab + ab\otimes 1 - b\otimes a
$$
satisfies the constant QYBE,
we seek a solution to equation~(\ref{yb}) of the following form
\begin{equation}\label{rans}
R(u,v)(a\otimes b) =\alpha(u,v)1\otimes ab + \beta(u,v)ab\otimes 1
-\gamma(u,v)b\otimes a,
\end{equation}
where $\alpha, \beta,\gamma$ are $k$-valued functions on $X\times X$.

Now, inserting this ansatz into equation~(\ref{yb}) one finds that
$R(u,v)$ is a solution of the two-parameter QYBE if and
only if
the coefficients of certain terms are equal to zero.
We present some of these calculations below.

\bigskip

$ R^{uv}_{12} \circ R^{uw}_{13} \circ R^{vw}_{23} (a \otimes b \otimes c) =
\beta(u,v)\beta(u,w)\beta(v,w) abc \otimes 1 \otimes 1 + $

$ \alpha(u,v)\beta(u,w)\beta(v,w) 1 \otimes abc \otimes 1 -
\gamma(u,v)\beta(u,w)\beta(v,w) bc \otimes a \otimes 1 + $

$ \beta(u,v)\alpha(u,w)\beta(v,w) bc \otimes 1 \otimes a +
\alpha(u,v)\alpha(u,w)\beta(v,w) 1 \otimes bc \otimes a - $

$ \gamma(u,v)\alpha(u,w)\beta(v,w) bc \otimes 1 \otimes a -
\beta(u,v)\gamma(u,w)\beta(v,w) bc \otimes 1 \otimes a - $

$ \alpha(u,v)\gamma(u,w)\beta(v,w) 1 \otimes bc \otimes a +
\gamma(u,v)\gamma(u,w)\beta(v,w) bc \otimes 1 \otimes a +
$

$
\beta(u,v)\beta(u,w)\alpha(v,w) abc \otimes 1 \otimes 1 +
\alpha(u,v)\beta(u,w)\alpha(v,w) 1 \otimes abc \otimes 1 -
$

$
\gamma(u,v)\beta(u,w)\alpha(v,w) 1 \otimes abc \otimes 1 +
\beta(u,v)\alpha(u,w)\alpha(v,w) 1 \otimes 1 \otimes abc +
$

$
\alpha(u,v)\alpha(u,w)\alpha(v,w) 1 \otimes 1 \otimes abc -
\gamma(u,v)\alpha(u,w)\alpha(v,w) 1 \otimes 1 \otimes abc -
$

$
\beta(u,v)\gamma(u,w)\alpha(v,w) bc \otimes 1 \otimes a -
\alpha(u,v)\gamma(u,w)\alpha(v,w) 1 \otimes bc \otimes a +
$

$
\gamma(u,v)\gamma(u,w)\alpha(v,w) 1 \otimes bc \otimes a -
\beta(u,v)\beta(u,w)\gamma(v,w) abc \otimes 1 \otimes 1 -
$

$
\alpha(u,v)\beta(u,w)\gamma(v,w) 1 \otimes abc \otimes 1 +
\gamma(u,v)\beta(u,w)\gamma(v,w) c \otimes ab \otimes 1 -
$

$
\beta(u,v)\alpha(u,w)\gamma(v,w) c \otimes 1 \otimes ab -
\alpha(u,v)\alpha(u,w)\gamma(v,w) 1 \otimes c \otimes ab +
$

$
\gamma(u,v)\alpha(u,w)\gamma(v,w) c \otimes 1 \otimes ab +
\beta(u,v)\gamma(u,w)\gamma(v,w) bc \otimes 1 \otimes a +
$

$
\alpha(u,v)\gamma(u,w)\gamma(v,w) 1 \otimes bc \otimes a -
\gamma(u,v)\gamma(u,w)\gamma(v,w) c \otimes b \otimes a  $

\bigskip

$ R^{vw}_{23} \circ R^{uw}_{13} \circ R^{uv}_{12} (a \otimes b \otimes c) =
\beta(v,w)\beta(u,w)\beta(u,v) abc \otimes 1 \otimes 1 + $

$
\alpha(v,w)\beta(u,w)\beta(u,v) abc \otimes 1 \otimes 1 -
\gamma(v,w)\beta(u,w)\beta(u,v) abc \otimes 1 \otimes 1
$

$
\beta(v,w)\alpha(u,w)\beta(u,v) 1 \otimes abc \otimes 1 +
\alpha(v,w)\alpha(u,w)\beta(u,v) 1 \otimes 1 \otimes abc -
$

$
\gamma(v,w)\alpha(u,w)\beta(u,v) 1 \otimes abc \otimes 1 -
\beta(v,w)\gamma(u,w)\beta(u,v) c \otimes ab \otimes 1 -
$

$
\alpha(v,w)\gamma(u,w)\beta(u,v) c \otimes 1 \otimes ab +
\gamma(v,w)\gamma(u,w)\beta(u,v) c \otimes ab \otimes 1
$

$
\beta(v,w)\beta(u,w)\alpha(u,v) c \otimes ab \otimes 1 +
\alpha(v,w)\beta(u,w)\alpha(u,v) c \otimes 1 \otimes ab -
$

$
\gamma(v,w)\beta(u,w)\alpha(u,v) c \otimes 1 \otimes ab +
\beta(v,w)\alpha(u,w)\alpha(u,v) 1 \otimes abc \otimes 1 +
$

$
\alpha(v,w)\alpha(u,w)\alpha(u,v) 1 \otimes 1 \otimes abc -
\gamma(v,w)\alpha(u,w)\alpha(u,v) 1 \otimes c \otimes ab -
$

$
\beta(v,w)\gamma(u,w)\alpha(u,v) c \otimes ab \otimes 1-
\alpha(v,w)\gamma(u,w)\alpha(u,v) c \otimes 1 \otimes ab +
$

$
\gamma(v,w)\gamma(u,w)\alpha(u,v) c \otimes 1 \otimes ab -
\beta(v,w)\beta(u,w)\gamma(u,v) bc \otimes a \otimes 1 -
$

$
\alpha(v,w)\beta(u,w)\gamma(u,v) bc \otimes 1 \otimes a +
\gamma(v,w)\beta(u,w)\gamma(u,v) bc \otimes 1 \otimes a -
$

$
\beta(v,w)\alpha(u,w)\gamma(u,v) 1 \otimes abc \otimes 1 -
\alpha(v,w)\alpha(u,w)\gamma(u,v) 1 \otimes 1 \otimes abc +
$

$
\gamma(v,w)\alpha(u,w)\gamma(u,v) 1 \otimes bc \otimes a +
\beta(v,w)\gamma(u,w)\gamma(u,v) c \otimes ab \otimes 1 +
$

$
\alpha(v,w)\gamma(u,w)\gamma(u,v) c \otimes 1 \otimes ab -
\gamma(v,w)\gamma(u,w)\gamma(u,v) c \otimes b \otimes a
$

\bigskip

The equality
$ R^{uv}_{12} \circ R^{uw}_{13} \circ R^{vw}_{23} (a \otimes b \otimes c) =
R^{vw}_{23} \circ R^{uw}_{13} \circ R^{uv}_{12} (a \otimes b \otimes c) $
implies that the coefficients of the following terms are equal to zero:
$1 \otimes abc \otimes 1 $,
$ bc \otimes 1 \otimes a $, $ 1 \otimes bc \otimes a $,
$ c \otimes ab \otimes 1 $ and $ c \otimes 1 \otimes ab $.
Thus, we obtain the following system of equations:
\begin{eqnarray} &&
(\beta(v,w)-\gamma(v,w))(\alpha(u,v)\beta(u,w) -
\alpha(u,w)\beta(u,v))\nonumber \\ &&\quad \quad \quad +
(\alpha(u,v)-\gamma(u,v))(\alpha(v,w)\beta(u,w) - \alpha(u,w)\beta(v,w))
= 0 \label{e1} \\ \nonumber \\ &&
\beta(v,w)(\beta(u,v)-\gamma(u,v))(\alpha(u,w)-\gamma(u,w)) \nonumber \\
&&\quad \quad \quad +
(\alpha(v,w)-\gamma(v,w))(\beta(u,w)\gamma(u,v)-\beta(u,v)\gamma(u,w)) =
0 \label{e2} \\ \nonumber \\ && \alpha(u,v)
\beta(v,w)(\alpha(u,w)-\gamma(u,w)) + \alpha(v,w)\gamma(u,w)
(\gamma(u,v) - \alpha(u,v)) \nonumber \\ &&\quad \quad \quad +
\gamma(v,w) (\alpha(u,v)\gamma(u,w)-\alpha(u,w)\gamma(u,v)) = 0
\label{e3} \\ \nonumber \\ && \alpha(u,v)
\beta(v,w)(\beta(u,w)-\gamma(u,w)) + \beta(v,w)\gamma(u,w) (\gamma(u,v)
- \beta(u,v)) \nonumber \\ &&\quad \quad \quad + \gamma(v,w)
(\beta(u,v)\gamma(u,w)-\beta(u,w)\gamma(u,v)) = 0 \label{e4} \\
\nonumber \\ && \alpha(u,v)(\alpha(v,w)-\gamma(v,w))(\beta(u,w) -
\gamma(u,w)) \nonumber \\ &&\quad \quad \quad + (\beta(u,v)-
\gamma(u,v))( \alpha(u,w) \gamma(v,w) - \alpha(v,w) \gamma(u,w)) = 0 \
\label{e5} \end{eqnarray}

The system of equations (\ref{e1}--\ref{e5}) is rather non-trivial.
Nonetheless,
it has some remarkable symmetry properties which can be used to
find some solutions. For example, let us observe that the equations
(\ref{e2}) and (\ref{e5}) are in some sense dual to each other.
Likewise, (\ref{e3}) and (\ref{e4}) are in some sense dual to each
other. In this paper we find some families of solutions for this
system. It is an open problem to classify all its solutions.

First, we try to find solutions to the
system of equations (\ref{e1}--\ref{e5})
of the following form:
$ \  \alpha(u,v)= pu-p'v $, $ \ \beta(u,v)=qu-q'v \ $ and
$ \gamma (u,v)= ru-r'v $.
We obtain the solutions:
$ \  \alpha(u,v)= p(u-v) $, $ \ \beta(u,v)=q(u-v) \ $ and
$ \gamma (u,v)= pu-qv $.
Thus, we are guided to the following results.

\begin{theorem} \label{thm1}
i) For any two parameters $p,q\in k$, the function\\
$R:X\times X\to \rend k {A\otimes A}$ defined by
\begin{equation}\label{rsol}
R(u,v)(a\otimes b) =p(u-v)1\otimes ab + q(u-v)ab\otimes 1 -(pu-qv)b\otimes a,
\end{equation}
is a coloured Yang-Baxter operator.

ii) If $ \ pu \neq qv $ and $ \ qu \neq pv $ then the operator
(\ref{rsol}) is invertible. Moreover, the following formula holds:
$$ R^{-1}(u,v)(a\otimes b) = \frac{p(u-v)}{(qu-pv)(pu-qv)}ba\otimes 1 +
\frac{q(u-v)}{(qu-pv)(pu-qv)}1\otimes ba - \frac{1}{(pu-qv)}b\otimes a
\ . $$
\end{theorem}
\begin{proof}
i) The proof can be done by following the steps presented above.
Another way to prove this theorem is by direct calculations.
The computations are quite involved.

ii) This can be verified by direct calculations.
\end{proof}

\begin{remark}
For any coalgebra $( C, \Delta, \varepsilon )$ and
two parameters $p,q\in k$, the function\\
$R_{C}:X\times X\to \rend k {A\otimes A}$ defined by
\begin{equation}\label{rsol1}
R_{C}(u,v)(c\otimes d) =p(u-v) \varepsilon(c) \Delta(d)
+ q(u-v) \varepsilon(d) \Delta(c) - (pu-qv)d\otimes c,
\end{equation}
is a coloured Yang-Baxter operator.
This follows from section 2, { \em The transfer of the theory to
coalgebras },
of \cite{DasNic:yan}.

\end{remark}

The system of equations (\ref{e1}--\ref{e5}) has another class
of solutions:
$ \alpha(u,v)= p^u q^v \ $, $ \ \beta(u,v)=p^u s^v \ $ and
$ \gamma (u,v)= p^u s^v $.
Thus, we obtain the following theorem.

\begin{theorem}\label{thm2}
i) For any three parameters $p,q,s\in k$, the function\\
$R':X\times X\to \rend k {A\otimes A}$ defined by
\begin{equation}\label{rsol2}
R'(u,v)(a\otimes b) = p^u (q^v 1\otimes ab + s^v ab\otimes 1 - s^v b\otimes a)
\end{equation}
is a coloured Yang-Baxter operator.

ii) If $ p \neq 0, \ q \neq 0, \ s \neq 0 $, then the operator
(\ref{rsol2}) is invertible. Moreover, the following formula holds:
$$ R'^{-1}(u,v)(a\otimes b) = \frac{1}{p^{u}}( \frac{1}{s^{v}} ba
\otimes 1 + \frac{1}{q^{v}} 1 \otimes ba - \frac{1}{s^{v}} b \otimes
a ) \ . $$
\end{theorem}
\begin{proof} Similar to that of Theorem \ref{thm1}.
\end{proof}

\begin{remark}
The system of equations (\ref{e1}--\ref{e5}) has also the following class
of solutions:
$ \alpha(u,v)= p^u q^v \ $, $ \ \beta(u,v)=s^u q^v \ $ and
$ \gamma (u,v)= p^u q^v $.
\end{remark}

\bigskip

\section{SOLUTIONS IN DIMENSIONS 2 AND 3}

\smallskip

We consider the algebra
$ A = \frac{ k[X]}{(X^2-\sigma)} $, where $ \sigma \in \{0, 1\}$ is a scalar.
Then $A$ has the basis $\{1, x \}$, where $x$ is the image of $X$ in the factor
ring. We consider the basis $ \{ 1 \otimes 1, 1 \otimes x, x \otimes 1,
x \otimes x \} $ of $ A \otimes A $ and represent the operator (\ref{rsol}) in
this basis:
\begin{eqnarray*}
&& R(u,v)(1\otimes 1) = (qu-pv) 1 \otimes 1 \\
&& R(u,v)(1\otimes x) = p(u-v) 1 \otimes x + (q-p)u x \otimes 1 \\
&& R(u,v)(x\otimes 1) = (q-p)v 1 \otimes x + q(u-v) x \otimes 1 \\
&& R(u,v)(x\otimes x) = \sigma (p+q)(u-v) 1 \otimes 1 -
(pu-qv) x \otimes x
\end{eqnarray*}
In matrix form, this operator reads
\begin{equation} \label{rmat}
R(u,v)= \begin{pmatrix}
qu-pv & 0 & 0 & \sigma (q+p)(u-v)\\
0 & p(u-v) & (q-p)v & 0\\
0 & (q-p)u & q(u-v) & 0\\
0 & 0 & 0 & qv-pu
\end{pmatrix}
\end{equation}
and satisfies the coloured QYBE. This equips us to look at the
FRT bialgebra
structure associated to this $R$--matrix operator. Let us recall
\cite{deeps,quesne} (and references therein) that the coloured extension
of a FRT bialgebra involves the generators to be parametrized by some
continuously varying {\em colour} parameters, and redefining the algebra
and the coalgebra such that all Hopf algebraic properties remain preserved.
Using the coloured FRT approach, the coloured $RTT$-- relations are
\begin{equation} \label{rtt}
R(u,v)T_{1u}T_{2v}=T_{2v}T_{1u}R(u,v)
\end{equation}
where $T_{1u}=T_{u}\otimes \bf{1}$, $T_{2v}={\bf 1} \otimes T_{v}$.
The generators are arranged in the matrices
$$
T_{u}=\begin{pmatrix}
a_{u} & b_{u}\\
c_{u} & d_{u}
\end{pmatrix}, \qquad
T_{v}=\begin{pmatrix}
a_{v} & b_{v}\\
c_{v} & d_{v}
\end{pmatrix}
$$
The $RTT$--equation (\ref{rtt}) then gives the algebra commutation
relations:
\begin{eqnarray}
&& p(u-v)[a_{u},d_{v}]-(q-p)(uc_{v}b_{u}-vc_{u}b_{v})=0 \label{r1} \\
&& p(u-v)a_{u}c_{v}-(qu-pv)c_{v}a_{u}+(q-p)vc_{u}a_{v}=0\label{r2} \\
&& (qu-pv)a_{u}b_{v}-p(u-v)b_{v}a_{u}-(q-p)ua_{v}b_{u}+\sigma
(q+p)(u-v)c_{u}d_{v}=0\label{r3} \\
&& p(u-v)b_{v}c_{u}-q(u-v)c_{u}b_{v}-(q-p)u(a_{u}d_{v}-a_{v}d_{u})=0
\label{r4} \\
&& (qv-pu)c_{u}d_{v}-p(u-v)d_{v}c_{u}-(q-p)uc_{v}d_{u}=0\label{r5} \\
&& p(u-v)b_{u}d_{v}-(qv-pu)d_{v}b_{u}+(q-p)vd_{u}b_{v}-\sigma
(q+p)(u-v)c_{v}a_{u}=0 \label{r6}
\end{eqnarray}
\begin{eqnarray}
&& (qu-pv)[a_{u},a_{v}]+\sigma (q+p)(u-v)c_{u}c_{v}=0\label{r7}
\\
&& (qu-pv)b_{u}b_{v}-(qv-pu)b_{v}b_{u}-\sigma
(q+p)(u-v)(a_{v}a_{u}-d_{u}d_{v})=0\label{r8} \\
&& (qv-pu)c_{u}c_{v}-(qu-pv)c_{v}c_{u}=0\label{r9} \\
&& [d_{u},d_{v}]=-[a_{u},a_{v}]\label{r10} \\
&& [a_{u},d_{v}]=[a_{v},d_{u}]\label{r11} \\
&& q(vc_{u}b_{v}-uc_{v}b_{u})-p(vb_{v}c_{u}-ub_{u}c_{v})=0 \label{r12}
\end{eqnarray}
All the above relations (\ref{r1}--\ref{r12}) are symmetric
with respect to the $u\leftrightarrow v$ exchange. Note that relation
(\ref{r9}) can be obtained from (\ref{r7}) when $\sigma \neq 0$,
but holds independently of (\ref{r7}) if $\sigma =0$.
This algebra does not have an uncoloured counterpart i.e. the limit $u=v$
does not hold. However, an interesting algebra arises in the limiting
case of $p=q$.
The coproduct is $\Delta(T_{u})=T_{u}\dot{\otimes}T_{u}$ where
$\dot{\otimes}$ is the usual symbol for tensor product and matrix
multiplication at the same time. The counit is $\varepsilon \left(
\begin{smallmatrix}a_{u}&b_{u}\\c_{u}&d_{u}\end{smallmatrix}\right)=
\left(\begin{smallmatrix}\bf{1}&0\\0&\bf{1}\end{smallmatrix}\right)$.
Finding a quantum determinant remains a difficult task, and we
conjecture that the above algebra does not admit one.

In the limit $p=q$, this algebra reduces to
\begin{equation*}
\begin{array}{lll}
& [a_{u},d_{v}]=0 \qquad & [a_{u},b_{v}]+2\sigma c_{u}d_{v}=0\\
& [a_{u},c_{v}]=0 \qquad & \{ b_{u},d_{v}\} -2\sigma c_{v}a_{u}=0\\
& [b_{v},c_{u}]=0 \qquad & [a_{u},a_{v}]+2\sigma c_{u}c_{v}=0\\
& \{ c_{u},d_{v}\} =0 \qquad & \{ b_{u},b_{v}\} -2\sigma
(a_{v}a_{u}-d_{u}d_{v})=0\\
& \{ c_{u},c_{v}\} =0 \qquad &
\end{array}
\end{equation*}
while relations (\ref{r10}) and (\ref{r11}) remain unchanged. Note that
the relation (\ref{r12}) arises as a compatibility condition for the
$u\leftrightarrow v$ exchange symmetry, which does not have an analogue
when $p=q$. This limit is interesting since the deformation parameter $q$
can be factored out from the $R$--matrix (\ref{rmat}) itself which then
depends only on the difference $u-v$. Consequently, the algebra in this
limit is also independent of $q$.

The operator (\ref{rsol2}) in
the same basis reads
\begin{equation} \label{rmat2}
R'(u,v)=  p^{u} \begin{pmatrix}
q^{v} & 0 & 0 & \sigma (q^{v}+ s^{v})\\
0 & q^{v} & q^{v}-s^{v} & 0\\
0 & 0 & s^{v} & 0\\
0 & 0 & 0 & -s^{v}
\end{pmatrix}
\end{equation}

Next, we present solutions in dimension three.
Consider the algebra
$ B = \frac{ k[X]}{(X^3-\eps X-\rho)} $, where $\eps $ and $ \rho $
are scalars.
Then $B$ has the basis $\{1, x, x^{2} \}$, where $x$ is the image of
$X$ in the factor
ring. We consider the basis $ \{ 1 \otimes 1, 1 \otimes x, 1 \otimes x^{2},
x \otimes 1, x \otimes x, x \otimes x^{2},
x^{2} \otimes 1, x^{2} \otimes x, x^{2} \otimes x^{2} \} $ of $ B \otimes
B $
and represent the operator (\ref{rsol}) in
this basis, in the matrix form
\begin{equation}\label{rmat3}
R(u,v)=\begin{pmatrix}
w & 0 & 0 & 0 & 0 & \rho t' \lambda & 0 & \rho t' \lambda & 0\\
0 & p \lambda & 0 & tv & 0 & \eps p \lambda & 0 & \eps p \lambda &
\rho p \lambda\\
0 & 0 & p \lambda & 0 & p \lambda & 0 & tv & 0 & \eps p \lambda \\
0 & tu & 0 & q \lambda & 0 & \eps q \lambda & 0 & \eps q \lambda & \rho
q \lambda\\
0 & 0 & 0 & 0 & w' & 0 & 0 & 0 & 0\\
0 & 0 & 0 & 0 & 0 & 0 & 0 & w' & 0\\
0 & 0 & tu & 0 & q \lambda & 0 & q \lambda & 0 & \eps q \lambda\\
0 & 0 & 0 & 0 & 0 & w' & 0 & 0 & 0\\
0 & 0 & 0 & 0 & 0 & 0 & 0 & 0 & w'
\end{pmatrix}
\end{equation}
where
$ \lambda = u-v, \  t=q-p, \ t'=q+p, \ w=qu-pv $, and $ w'=qv-pu $.

Similarly, the operator (\ref{rsol2}) can be represented in the same
basis and has the matrix form
\begin{equation}
R'(u,v)= p^{u} \begin{pmatrix}
q^{v} & 0 & 0 & 0 & 0 & \rho (q^{v} + s^{v}) & 0 & \rho (q^{v} + s^{v}) & 0\\
0 & q^{v} & 0 & q^{v} - s^{v} & 0 & \eps q^{v} & 0 & \eps q^{v} & \rho q^{v}\\
0 & 0 & q^{v} & 0 & q^{v} & 0 & q^{v} - s^{v} & 0 & \eps q^{v}\\
0 & 0 & 0 & s^{v} & 0 & \eps s^{v} & 0 & \eps s^{v} & \rho s^{v}\\
0 & 0 & 0 & 0 & -s^{v} & 0 & 0 & 0 & 0\\
0 & 0 & 0 & 0 & 0 & 0 & 0 & -s^{v} & 0\\
0 & 0 & 0 & 0 & s^{v} & 0 & s^{v} & 0 & \eps s^{v}\\
0 & 0 & 0 & 0 & 0 & -s^{v} & 0 & 0 & 0\\
0 & 0 & 0 & 0 & 0 & 0 & 0 & 0 & -s^{v}
\end{pmatrix}
\end{equation}

%\smallskip

\section{One--parameter Yang--Baxter operators}

\smallskip

The analysis presented above leads naturally to look at solutions of the
one-parameter QYBE. Let $X$ be a set, $Z \subset X \times X$ and $V$
is a finite dimensional vector space over a field $k$.
Formally, an one-parameter Yang-Baxter operator is defined as a function $$ R
:X \to \rend k {V\otimes_k V}. $$ Thus for any $x\in X$,
$R(x) : V\otimes_k V\to V\otimes_k V$ is a linear operator.
$R$ is a one-parameter Yang-Baxter operator if it
satisfies the one-parameter form of the QYBE,
\begin{equation}\label{oneparyb}
R_{12}(x)R_{13}(\varphi (x,z))R_{23}(z) = R_{23}(z)
R_{13}(\varphi (x,z))R_{12}(x)
\end{equation}
for all $(x,z) \in Z$, where $ \varphi : Z \to X $.

We seek a
solution to equation~(\ref{oneparyb}) of the form
\begin{equation}\label{rans}
R(x)(a\otimes b) =\alpha(x)1\otimes ab + \beta(x)ab\otimes 1
-\gamma(x)b\otimes a,
\end{equation}
where $\alpha, \beta,\gamma$ are $k$-valued functions on $X$.

Now, inserting this ansatz into equation~(\ref{oneparyb}) one finds
that $R(x)$ is a solution of the one-parameter QYBE if and only if
the coefficients of certain terms are equal to zero. Thus, we obtain
the following system of equations:

\begin{eqnarray} &&
(\beta(z)-\gamma(z))(\alpha(x)\beta(\varphi (x,z)) -
\alpha(\varphi (x,z))\beta(x))\nonumber \\ &&\quad \quad \quad +
(\alpha(x)-\gamma(x))(\alpha(z)\beta(\varphi (x,z)) - \alpha(\varphi
(x,z))\beta(z))
= 0 \label{f1} \\ \nonumber \\ &&
\beta(z)(\beta(x)-\gamma(x))(\alpha(\varphi (x,z))-\gamma(\varphi
(x,z))) \nonumber \\
&&\quad \quad \quad +
(\alpha(z)-\gamma(z))(\beta(\varphi
(x,z))\gamma(x)-\beta(x)\gamma(\varphi (x,z))) =
0 \label{f2} \\ \nonumber \\ &&
\alpha(x) \beta(z)(\alpha(\varphi
(x,z))-\gamma(\varphi (x,z))) + \alpha(z)\gamma(\varphi (x,z))
(\gamma(x) - \alpha(x)) \nonumber \\ &&\quad \quad \quad +
\gamma(z) (\alpha(x)\gamma(\varphi (x,z))-\alpha(\varphi (x,z))\gamma(x)) = 0
\label{f3} \\ \nonumber \\ &&
\alpha(x) \beta(z)(\beta(\varphi
(x,z))-\gamma(\varphi (x,z))) + \beta(z)\gamma(\varphi (x,z))
(\gamma(x)
- \beta(x)) \nonumber \\ &&\quad \quad \quad
+ \gamma(z)
(\beta(x)\gamma(\varphi (x,z))-\beta(\varphi (x,z))\gamma(x)) = 0
\label{f4} \\
\nonumber \\ && \alpha(x)(\alpha(z)-\gamma(z))(\beta(\varphi (x,z)) -
\gamma(\varphi (x,z))) \nonumber \\ &&\quad \quad \quad + (\beta(x)-
\gamma(x))( \alpha(\varphi (x,z)) \gamma(z) - \alpha(z)
\gamma(\varphi (x,z))) = 0 \
\label{f5} \end{eqnarray}

Finding further solutions for this system of
equations and their classification is also an open problem.  From the above
analysis it follows that the
system of equations (\ref{f1}--\ref{f5})
has the following solution:
$ \  \alpha(x)= x-1 $, $ \ \beta(x)=q(x-1) \ $, and
$ \gamma (x)= x-q $, where $ \varphi (x,z)=xz $.
Thus, we are guided to the following results.

\begin{proposition} \label{prop1}
i) For any parameter $q\in k$, the function $R:X \to \rend k
{A\otimes A}$ defined by
\begin{equation}\label{oneparsol}
R(x)(a\otimes b) =(x-1)1\otimes ab + q(x-1)ab\otimes 1 -(x-q)b\otimes a,
\end{equation}
is a one-parameter Yang-Baxter operator, where $ \varphi (x,z)=xz $.

ii) If $ \ x \neq q $ and $ \ qx \neq 1 $, then the operator (\ref{oneparsol})
is invertible. Moreover, the following formula holds:
$$ R^{-1}(x)(a\otimes b) = \frac{x-1}{(qx-1)(x-q)}ba\otimes 1 +
\frac{q(x-1)}{(qx-1)(x-q)}1\otimes ba - \frac{1}{(x-q)}b\otimes a
\ . $$

iii) For any coalgebra $( C, \Delta, \varepsilon )$ and
a parameter $q\in k$, the function
$R_{C}:X \to \rend k {C\otimes C}$ defined by
\begin{equation}\label{rsol1}
R_{C}(x)(c\otimes d) =(x-1) \varepsilon(c) \Delta(d)
+ q(x-1) \varepsilon(d) \Delta(c) - (x-q)d\otimes c,
\end{equation}
is a one-parameter Yang-Baxter operator, where $ \varphi (x,z)=xz $.
\end{proposition}

\begin{remark}
We consider the algebra $ A = \frac{
k[X]}{(X^2-\sigma)} $ as before and the same basis of $A$.
Then the
operator  (\ref{oneparsol}) reads

\begin{equation} \label{oneparrmat}
R(x)= \begin{pmatrix}
qx-1 & 0 & 0 & \sigma (q+1)(x-1)\\
0 & (x-1) & (q-1) & 0\\
0 & (q-1)x & q(x-1) & 0\\
0 & 0 & 0 & q-x
\end{pmatrix}
\end{equation}
\end{remark}

The system of equations (\ref{f1}--\ref{f5}) has another class
of solutions:
$ \alpha(x)= x  $, $ \ \beta(x)= 1 \ $ and
$ \gamma (x)= 1 $, where $ \varphi (x,z)= z $.
Thus, we obtain the following proposition.

\begin{proposition}\label{proponepar}
i) The function $R':X \to \rend k {A\otimes A}$ defined by
\begin{equation}\label{oneparrsol2}
R'(x)(a\otimes b) = x 1\otimes ab + ab\otimes 1 - b\otimes a
\end{equation}
is a one-parameter Yang-Baxter operator, where $ \varphi (x,z)= z $.

ii) If $ x \neq 0 $, then the operator
(\ref{oneparrsol2}) is invertible. Moreover, the following formula holds:
$$ R'^{-1}(x)(a\otimes b) =   ba \otimes 1 + \frac{1}{x} 1 \otimes ba
- b \otimes a  \ . $$
\end{proposition}

\begin{remark}\label{remonepar}
The system of equations (\ref{f1}--\ref{f5}) has also the following class
of solutions:
$ \alpha(x)= 1$, $ \ \beta(x)= x \ $, and
$ \gamma (x)= 1 $, where $ \varphi (x,z)= x $.
\end{remark}

Let us now compare our solutions with other well-known spectral-parameter dependent solutions of the QYBE. Of particular interest is a solution appearing in the context of generalisation of some exactly solvable $q$--state vertex models  \cite{perk}. This solution coincides with another solution associated to a multiparametric quantum deformation of the universal enveloping algebra of a symmetrizable Kac-Moody algebra \cite{okado}. Consider the formula (3.3) in \cite{okado} (which is the multiparameter $R$--matrix for $U_{q,Q}(\widehat{sl}(n))$) for $n=2$. We obtain

\begin{eqnarray}
{\hat R}  =  && c_{11}(q^2x-1)E_{11}\otimes E_{11} + c_{22}(q^2x-1)E_{22}\otimes E_{22} \nonumber\\
&& + c_{22}\gamma_{12}(q^2-1)E_{11}\otimes E_{22} +  c_{11}\gamma_{21}(q^2-1)xE_{22}\otimes E_{11} \label{oyformula} \\
&& + c_{12}s_{12}q(x-1)E_{12}\otimes E_{21} +  c_{21}s_{21}q(x-1)E_{21}\otimes E_{12} \nonumber
\end{eqnarray}

where parameters $c$, $s$ and $\gamma$'s are defined in \cite{okado}. In line with Remark $2$ of \cite{okado}, we specialise these to be $1$ to obtain a single-parameter version. Note that (\ref{oyformula}) is a solution of the {\em braid} equation, which in the matrix form reads

\begin{equation} \label{okadormat}
{\hat R}= \begin{pmatrix}
q^2x-1 & 0 & 0 & 0\\
0 & q^2-1 & q(x-1) & 0\\
0 & q(x-1) & (q^2-1)x & 0\\
0 & 0 & 0 & q^2x-1
\end{pmatrix}
\end{equation}

On the other hand, we also note that our $R$--matrix (\ref{rmat}) reduces to (\ref{oneparrmat}) by setting $p=1$ and adjusting an overall factor of $\frac{1}{v}$. In order to compare (\ref{oneparrmat}) with (\ref{okadormat}) we have to set $\sigma=0$ and compose ({\ref{oneparrmat}) with the {\em twist} map to obtain it as a solution of the braid equation, from that of the QYBE (\ref{yb}). It turns out that both solutions (\ref{okadormat}) and twisted (\ref{oneparrmat}), though distinct,  are very similar to each other and if we further set $q=1$, then both reduce to the same trivial $R$--matrix. It would be interesting to find further classes of solutions to the system of equations (\ref{f1}--\ref{f5}) that might relate to other known solutions.

\smallskip

\section{Yang--Baxter systems}

\smallskip

Yang-Baxter systems were
introduced in \cite{HlaSno:sol} as a spectral-parameter independent
generalisation of the quantum Yang-Baxter equation related to
non-ultralocal integrable systems studied previously in
\cite{HlaKun:qua}. Yang-Baxter systems are conveniently defined in
terms of {\em Yang-Baxter commutators}. Consider three vector spaces
$V,V',V''$ and  three linear maps
$ R : V \ot V' \rightarrow V \ot V' $,
$ S : V \ot V'' \rightarrow V \ot V'' $, and
$ T : V' \ot V'' \rightarrow V' \ot V'' $. Then a {\em Yang-Baxter
commutator} is a map
$ [R,S,T]:  V \ot V' \ot V'' \rightarrow V \ot V' \ot V'' $, defined by
\begin{equation}   \label{ybcomm}
[R,S,T]= R_{12}  \circ  S_{13}  \circ  T_{23} - T_{23}  \circ  S_{13}
\circ  R_{12} \ .
\end{equation}
In terms of a Yang-Baxter commutator,  the quantum Yang-Baxter
equation (\ref{constyb})  is expressed simply as $[R,R,R] = 0$.

\begin{definition}\label{def.wxz}
Let $V$ and  $V'$ be vector spaces. A system of linear maps
$$ W : V \ot V \rightarrow V \ot V ,\quad
Z : V' \ot V' \rightarrow V' \ot V' , \quad
  X : V \ot V' \rightarrow V \ot V'
$$
is called  a {\em WXZ-system} or a
{\em Yang-Baxter system}, provided the following equations are satisified:
\begin{equation}   \label{ybeqn4}
[W,W,W]\ = \  0 \ ,
\end{equation}
\begin{equation}   \label{ybeqn5}
[Z,Z,Z]\ = \  0 \ ,
\end{equation}
\begin{equation}   \label{ybeqn6}
[W,X,X]\ = \  0 \ ,
\end{equation}
\begin{equation}   \label{ybeqn7}
[X,X,Z]\ = \  0 \ .
\end{equation}
\end{definition}

There are several algebraic origins and applications of WXZ-systems. It has been
observed in \cite{Vla:met} that WXZ-systems with invertible $W$, $X$ and $Z$
can be used to construct dually-paired bialgebras of the FRT type, thus leading
to quantum doubles. Given a WXZ-system as in Definition~\ref{def.wxz} one can
construct an invertible solution to the constant Yang-Baxter equation provided
W, X and Z are invertible (see \cite{BrzNic:yan} for details). In
\cite{BrzNic:yan} it was also shown that a Yang-Baxter system can be
constructed from any entwining structure and conversely, that
Yang-Baxter systems of certain types lead to entwining structures.

The next result is a new construction of a Yang-Baxter system.

\begin{theorem}\label{thm3}
Let $A$ be a k-algebra and $
\lambda, \mu \in k$.
The following is a Yang-Baxter system:

$ W : A \ot A \rightarrow A \ot A ,\quad
W (a\otimes b) = \lambda 1\otimes ab + ab\otimes 1 - b\otimes a $

$ Z : A \ot A \rightarrow A \ot A , \quad
Z (a\otimes b) =  1\otimes ab +  \mu ab\otimes 1 - b\otimes a $

$ X : A \ot A \rightarrow A \ot A, \quad
X (a\otimes b) =  1\otimes ab + ab\otimes 1 - b\otimes a $
\end{theorem}
\begin{proof}
This follows from Proposition \ref{proponepar},
Remark \ref{remonepar}, and \cite{DasNic:yan}.
\end{proof}

We now represent the Yang-Baxter
system from the theorem above in dimension two in the same basis as
in the previous sections:
\begin{equation} \label{ybst}
W= \begin{pmatrix}
\lambda & 0 & 0 & \sigma ( \lambda +1)\\
0 & \lambda &  \lambda -1 & 0\\
0 & 0 & 1 & 0\\
0 & 0 & 0 & -1
\end{pmatrix}
\end{equation}
\begin{equation} \label{ybst2}
Z= \begin{pmatrix}
\mu & 0 & 0 & \sigma ( \mu +1)\\
0 & 1 &  0 & 0\\
0 & \mu -1 & \mu & 0\\
0 & 0 & 0 & -1
\end{pmatrix}
\end{equation}
\begin{equation} \label{ybst3}
X= \begin{pmatrix}
1 & 0 & 0 & 2 \sigma\\
0 & 1 &  0 & 0\\
0 & 0 & 1 & 0\\
0 & 0 & 0 & -1
\end{pmatrix}
\end{equation}

\smallskip

\section{Concluding remarks}

\smallskip

One of our main results is Theorem \ref{thm1}, defining a class of coloured
Yang-Baxter operators (\ref{rsol}) and both (\ref{rmat}) and
(\ref{rmat3}) are the respective matrix forms, {\em i.e.}
spin--$\frac{1}{2}$ and spin--$1$ $R$--matrices, for most general algebras
in dimensions $2$ and $3$, respectively. This is distinct from the
solutions appeared in \cite{HlaWan}. (\ref{rmat}) is likely to correspond
to some seven-vertex solution of Baxter's solvable models. Similarly,
Theorem \ref{thm2} produces another class of such operators. As remarked
earlier, it is an open problem to classify all coloured Yang-Baxter
operators solving the system of equations (\ref{e1}--\ref{e5}).

It is pertinent to note that this work is in the spirit of references
\cite{DasNic:yan,Nic:sel}, {\em i.e.} seeking a coloured generalisation of
the (constant) Yang-Baxter operators that were derived from algebra
structures which
also produced new classes of solutions to the constant QYBE, again
distinct from those that had appeared in the literature (see for instance
\cite{Hiet}). This makes the framework of Yang-Baxter operators somewhat
different and perhaps more general than the traditional approach of
solving the QYBE.

We believe it could be of physical interest to look at Yang-Baxter
operators and their coloured generalisations to make contact with
integrability and vertex models in a way analogous to the additive
solutions of the spectral parameter dependent QYBE. Furthermore, we have
shown that one can associate a FRT bialgebra to a coloured Yang-Baxter
operator exhibiting explicitly for the $2$--dimensional case. It
would be useful to investigate further such algebras (including the
coloured quantum groups) which on the one hand seem quite constrained, but
on the other, possess nice properties of the colour exchange symmetries.

Proposition \ref{prop1} defines a class of one-parameter Yang-Baxter operators and
Theorem \ref{thm3} is another main result which associates a Yang-Baxter system
to such operators. Further work along these lines will lead to a
better understanding of spectral-parameter dependent Yang-Baxter
operators.

\section*{Acknowledgements}
We are grateful to Tomasz Brzezi\'nski for several fruitful
discussions. F. N. thanks the European Commission for the Marie
Curie Fellowship HPMF-CT-2002-01782 at the University of Wales Swansea,
and D. P. thanks the Royal Commission for the Exhibition of 1851 for a Research
Fellowship. D. P.'s research has been partly supported by a grant
from the European Science Foundation's programme
on Noncommutative Geometry. Finally, we thank the referee for drawing
our attention to references \cite{okado, perk}.

\end{document}